\documentclass[11pt]{article}
\usepackage[margin=1in, paperwidth=8.5in, paperheight=11in]{geometry}
\usepackage{graphicx}
\usepackage{float}
\usepackage{caption}
\usepackage{longtable}
\usepackage{booktabs} 
\usepackage{color}
\usepackage{amsmath}
\usepackage{cancel}
\usepackage{tikz}
\usepackage{amsfonts}
\usepackage[linesnumbered,ruled]{algorithm2e}

\usepackage{dsfont}

\usepackage[T1]{fontenc}
\usepackage[utf8x]{inputenc} 
\usepackage[affil-it]{authblk}
\usepackage{lmodern}
\usepackage[english]{babel}
\usepackage{amsfonts}
\usepackage{amssymb}
\usepackage{ stmaryrd }
\usepackage{ textcomp }
\usepackage{ bbold }
\usepackage{authblk}
\usepackage[toc,page]{appendix}
\usepackage{mathrsfs}
\usepackage[all]{xy} 
\usepackage{xcolor}
\usepackage{tikz} 
\usepackage{babel}
\usepackage{enumerate}  
\usepackage{amsthm}
\usepackage{verbatim}

\usepackage{setspace}
\usepackage{listings}
\PrerenderUnicode{é}
\usepackage{amsmath}

\newtheorem{theorem}{Theorem}[section]

\newtheorem{proposition}[theorem]{Proposition}

\providecommand{\keywords}[1]{\textbf{\textit{Keywords---}} #1}

\title{Optimal Transport Filtering with Particle Reweighing in Finance}
\author[1]{Raphael Douady}
\author[2]{Shohruh Miryusupov  \thanks{Corresponding author: \texttt{shohruh.miryusupov@malix.univ-paris1.fr}}}
\affil[1]{Stony Brook University, CNRS, Université Paris 1 Panthéon Sorbonne}
\affil[2]{Université Paris 1 Panthéon Sorbonne\\ Labex RéFi 
}

\begin{document}
\maketitle
\begin{abstract}
We show the application of an optimal transportation approach to estimate stochastic volatility process by using the flow that optimally transports the set of particles from the prior to a posterior distribution. We also show how to direct the flow to a rarely visited areas of the state space by using a particle method (a mutation and a reweighing mechanism). We demonstrate the efficiency of our approach on a simple example of the European option price under the Stein-Stein stochastic volatility model for which a closed form formula is available. Both homotopy and reweighted homotopy methods show a lower variance, root-mean squared errors and a bias compared to other filtering schemes recently developed in the signal-processing literature, including particle filter techniques. 
\end{abstract}
\keywords{optimal transport, Monge-Kantorovich, stochastic volatility, Particle methods, option pricing, Stein model, importance sampling, variance reduction, particle filter, Monte Carlo simulations, sequential Monte Carlo}

\tableofcontents

\newpage
\section{Introduction}

\hspace{3ex} Optimal transport problem that was formulated by Monge in the XVIII century, then reformulated fifty years ago by Kantorovich and it has recently been rediscovered by C. Villani \cite{CV}. This problem then was applied in different contexts, for example in option pricing \cite{HPL2013}. 

Particle methods which were extensively researched by P. Del Moral in \cite{DM2004}, \cite{DM2013} and \cite{DM} allow to find so-called "optimal transport". For this purpose a set of discrete  weighted samples, i.e. particles, is used to approximate an importance measure, and then to predict a posterior distribution by propagating the set of particles until we get an estimate. 

Another approach has been proposed by Daum's et al. \cite{DH2013}, \cite{DH2011} that allows the reduction of the number of particles we need in order to get a tolerable level of errors in the filtering problem. The main idea behind this method is the evolution in homotopy parameter $\lambda$ (a "pseudotime") from prior to the target density. They introduced a particle flow, in which particles are gradually transported without the necessity to randomly sample from any distribution. This approach as an optimal transport problem allows optimally move the set of particles according to Bayes' rule. In other words, the particles are progressively transported according to their flow. One can in this way reduce the number of needed samples, since the variance and bias of the estimator is lower and as a result reduce the computational burden in both the estimation and the prediction steps.

In this paper we adapt homotopy transport in Stein-Stein stochastic volatility model \cite{SS1991} to price a European option and extend Daum's et al. method by reweighing the generated particles' trajectories that allows to efficiently transport the particles from a prior transition density to a posterior one under the measurement impact. The idea of transportation and reweighing mechanism is to transport particles through the sequence of densities that move the least during the synthetic time until they reach the posterior distribution. By regenerating particles according to their weight at each time step we are able to direct the flow and further minimize the variance of the estimates. The transportation of particles can be understood as a geodesic flow in a convex subset of a Euclidean space.

We show that homotopy transport allows to significantly reduce the variance compared to a particle filtering technique. Path reweighing allows further reduce both the variance and the bias of estimators.

The rest of the article is organized as follows. Section 2 formulates the problem of computing the expectation  when we have partially observed variables and shows the solution using particle filter method. Section 3 formulates the problem defined in section 2 in the context of optimal transport and presents the homotopy transport approach to solve the problem. Section 4 shows the mixture of homotopy tranport and path reweighing approach and, actually, extends the method proposed in section 3. Section 5 provides numerical results. Section 6 concludes.

\section{Particle Filtering}
\subsection{Problem formulation}

\hspace{3ex} Many problems arises in financial applications when one has to compute expectations with partially observed information. A simple example is an option pricing with hidden volatility dynamics. Assume that we denote by $\{Y_t\}_{t\geq 0} \in \mathbb{R}^{n_Y}$ asset returns, which are observed from the dynamics of prices, while the hidden factor $\{X_t\}_{t\geq 0} \in \mathbb{R}^{n_X}$ is unobservable. Let $(\Omega, \mathcal{F}, \mathbb{P})$ be a probability space and the set of observed data at time $t$ be $(\mathcal{F}_t)$ is a filtration generated by a process $(Y_t)_{t\geq 0}$. 




The classical problem, where particle filtering is applied, is to extract a sequence of hidden variables $X_t$. It is formalized in the following way, given an initial $\mathbb{R}^{n_X}$-dimensional random variables $x_0$ with distribution $\mathbb{P}_{x_0}$, then for $t\in\mathbb{N}$:
\begin{equation}
\left\{ \begin{array}{c}
X_t = f(X_{t-1},\epsilon_t)\\
Y_t = h(X_t,Y_{t-1},\eta_t)
\end{array} \right.
\end{equation}
where the first equation is the hidden process, with  
$\epsilon_t:\Omega \rightarrow \mathbb{R}^{n_X}$ are i.i.d random variables, the map $f:\mathbb{R}^{n_X} \rightarrow \mathbb{R}^{n_X}$ is $\mathcal{B}(\mathbb{R}^{n_X})$ - measurable. The second equation is called a measurement model with $\eta_t:\Omega \rightarrow \mathbb{R}^{n_Y}$ are i.i.d. random variables and the map $h:\mathbb{R}^{n_X}\times\mathbb{R}^{n_Y} \rightarrow \mathbb{R}^{n_Y}$ is $\mathcal{B}(\mathbb{R}^{n_X})\otimes \mathcal{B}(\mathbb{R}^{n_Y})$ - measurable.  

Given above stochastic dynamical system, we would like to compute the following conditional expectation:
\begin{equation} \label{exp1}
\mathbb{E}[z(X_{t})|\mathcal{F}_t]=\frac{1}{\mathcal{Z}}\int \nu(dx_{0:t})\rho_t(x_{0:t},Y_{1:t})z(x_{t})
\end{equation}
with a distribution of $X_{0:t}$:
\begin{equation}
\nu(dx_{0:t})=p_0(dx_0)\prod_{l=0}^t k_l(X_{l-1},X_l)\mu(dx_l)
\end{equation}
and normalizing constant $\mathcal{Z}$:
\begin{equation}
\mathcal{Z}=\int \nu(dx_{0:t})\rho_t(x_{0:t},Y_{1:t})
\end{equation}
where $(X_t)_{t\geq 0}$ forms a Markov Chain in $(\mathbb{R}^{n_X},\mathcal{B}(\mathbb{R}^{n_X}))$ with transition density $k_t:\mathbb{R}^{n_X}\times \mathbb{R}^{n_X} \rightarrow \mathbb{R}^{n_X}_+$ with respect to the measure $\mu(dx)$. The random variables $(Y_t)_{t\geq 0}$ in $(\mathbb{R}^{n_Y},\mathcal{B}(\mathbb{R}^{n_Y}))$ are conditionally independent given $(X_t)_{t\geq 0}$ with transition density (likelihood) $\rho_t:\mathbb{R}^{n_X}\times \mathbb{R}^{n_Y} \rightarrow \mathbb{R}_+^{n_Y}$ with reference measure $\gamma$.



Intuitively, one could think that we could use naive Monte Carlo technique to approximate (\ref{exp1}):
\begin{equation}
\mathbb{E}[z(X_{t})|\mathcal{F}_t] \approx \int(\mathbb{M}^{n_X}\mathbb{P})(d (x_{0:t})z(x_t)=\frac{1}{n_X}\sum_{i=1}^{n_X} f(X^{(i)}_t)
\end{equation}
where the sampling operator $\mathbb{M}^{n_X}\nu=\frac{1}{n_X}\sum_{i=1}^{n_X}\delta_{X^{(i)}}$, $X \in \mathbb{R}^{n_X}$ and $\forall i=1,...,n_X$, $X^{(i)}$ are i.i.d. draws from $\nu$.

The problem with naive Monte Carlo Sampling lies in the fact that we don't know how to sample from conditional distribution $\mathbb{P}(d(x_{0:t}))=\mathbb{P}(X_{0:t}\in dx_{0:t}|Y_{1:t})$. Moreover, computation of normalization constant $\mathcal{Z}$ is a big challenge.

There is a lot of research made to tackle this problem, for example \cite{DM2004}, where the problem is transformed from a partially observed to a fully observed, by introducing a so called filtering distribution, that links observed and latent variables and recursively updates it. 


\begin{proposition}
Conditional probability (filtering distribution) $\Xi_t=P(X_t\in \cdot|Y_1,...,Y_t)$ with prior $X_0\sim p_0$ could be computed sequentially:
\begin{equation} \label{op1}
\Xi_tz=\frac{\int \Xi_{t-1}(dx_{t-1})k_t(X_{t-1},x)\mu(dx)\rho_t(x,Y_t)z(x)}{\int \Xi_{t-1}^{\mu}(dx_{t-1})k_t(X_{t-1},x)\mu(dx)\rho_t(x,Y_t)}
\end{equation}
with $\Xi_0=p_0$ and $\Xi_tz=\int P(x_{0:t}\in dx_{0:t}|Y_1,...,Y_t)z(x_{0:t})$
\end{proposition}

Denote the corresponding values of the hidden process as $(X_0,...,X_t)$  and the values of the measurement process as $(Y_0,...,Y_t)$. If there exists an absolutely continuous probability measure $\mathbb{P}\ll\mathbb{Q}$,  than for $t=0,...,N$ we have:
\begin{equation} \label{RDN1}
\mathbb{E}^{\mathbb{P}}[z(X_{t})|\mathcal{F}_t]=\mathbb{E}^{\mathbb{Q}}[z(X_{t})\frac{d\mathbb{P}}{d\mathbb{Q}}(X_{0:t})|\mathcal{F}_t]
\end{equation}

An importance measure $\mathbb{Q}$ could be chosen arbitrarily as soon as the continuity of the measure is preserved. But usually in a sequential importance sampling literature it is common to see the approximation of $\mathbb{Q}$, given that there exists an absolutely continuous importance kernel $\widetilde{K}_t$, such that for $K\ll \widetilde{K}_t$ as:
\begin{equation}
\mathbb{Q}(B)=\sum_{i=1}^M \omega^{(i)}_t \widetilde{K}_t(X^{(i)}_{t-1},A_i), \ \ B\in\mathcal{B}(\mathbb{R}^{n_X})
\end{equation}

where $A_i=\{X_{t}\in\mathbb{R}^{n_X}|\mathbb{1}_{B}(X^{(i)}_{t-1},X_{t})=1\}$, $(\omega_{t-1}^{(i)})_{i=1}^M$ is the weight function, and for $i=1,...,n_X$, $(X^{(i)}_0,...,X^{(i)}_t)$ are independent trajectory realizations. Now assume that the prior and sampling kernels $K_t$ and $\widetilde{K}_t$ have densities $k_t$ and $\widetilde{k}_t$ with respect to the measure $\mu$,  $\forall t=1,...,T$. 

For $0<...<t$, the Radon-Nikodym derivative in (\ref{RDN1}) is: 
\begin{equation}
\frac{d\mathbb{P}}{d\mathbb{Q}}(X_{0:t})=\frac{1}{\mathcal{Z}}\rho_1(X_1,Y_1)\frac{k_1(X_0,X_1)}{\widetilde{k}_1(X_0,X_1)}...\rho_t(X_{t},Y_t)\frac{k_t(X_{t-1},X_t)}{\widetilde{k}_{t}(X_{t-1},X_t)}
\end{equation}

where the importance measure is given by:
\begin{equation}
\mathbb{Q}(dx_{0:t})=p_0(dx_0)\widetilde{k}_1(x_0,x_1)\mu(dx_1)...\widetilde{k}_t(x_{t-1},x_t)\mu(dx_t)
\end{equation}

Observe, that we still can not compute a normalization constant $Z$, otherwise to compute the filtering distribution $\Xi_t$ will not be a problem, so we will need to apply normalized operator $\mathbb{M}^{N_x}$ to approximate filtering distribution:

\begin{equation}
\mathbb{E}^{\mathbb{P}}[z(X_{0:t})|\mathcal{F}_t]\approx \int\mathbb{M}^{n_X}\mathbb{P}(dx_{0:t})z(x_t)=\sum_{i=1}^{n_X}\widetilde{\omega}^{(i)}_tz(X_t^{(i)})\delta_{X_t^{(i)}}(dx_t)
\end{equation}

where the normalized importance weight function:
\begin{equation} \label{w11}
\widehat{\omega}_{t}^{(i)}(X_{t}^{(i)})=\frac{\omega_{t}^{(i)}(X_{t}^{(i)})}{\sum_{i=1}^M \omega_{t}^{(j)}(X_{t}^{(j)})}
\end{equation}

and an unnormalized weight is given by:

\begin{equation}
\omega_{t}^{(i)}(X_{t}^{(i)})=\prod_{l=1}^{t}\rho_{l}(X_{l}^{(i)},Y_{l})\frac{ k_l(X_{t-1}^{(i)},X_{l}^{(i)})}{\widetilde{k}_l(X_{l-1}^{(i)},X_{l}^{(i)})}
\end{equation} 
Observe that importance weights $\{\widehat{\omega}_{t}^{(i)}\}_{i=1}^{n_X}$ are positive and $\sum_{i=1}^{n_X}\widehat{\omega}_{t}^{(i)}=1$.

Since Particle filters showed weight degeneracey as number of time steps increased, Gordon et al. (1993) proposed a resampling step to the algorithm, which could be described by the following nonlinear equation:
\begin{equation}
\Xi_t=\Phi_t\Xi_{t-1}^{\mu} \ \ \mbox{with} \ \ \Xi_{0}=p_0
\end{equation}
where the nonlinear operator $\Phi_t$ is given by:
\begin{equation} 
(\Phi_t\nu)z=\frac{\int \nu(dx_{t-1})k_t(X_{t-1},x)\mu(dx)\rho_t(x,Y_t)z(x)}{\int \nu(dx_{t-1})k_t(X_{t-1},x)\mu(dx)\rho_t(x,Y_t)}
\end{equation}

The action of the operator $\Phi$  could be schematically described as:

\begin{equation}
\Xi_{t-1}  \xrightarrow[]{Mutation} \mathcal{M}\Xi_{t-1}  \xrightarrow{Reweighing} \Omega_t \mathcal{M}\Xi_{t-1}
\end{equation}

where the mutation operator $\mathcal{M}$ is given by

\begin{equation}
(\mathcal{M}\nu)(z)=\int \nu(dx_{t-1})p(x_{t-1},x)\mu(dx)z(x)
\end{equation}

and the reweighing operator $\Omega_t$ has the form

\begin{equation}
\Omega_t(\nu)z=\frac{\int \nu(dx)\rho_t(x,Y_t)f(x)}{\int \nu(dx)g(x,Y_t)}
\end{equation}

After the reweighing step we get the following approximation of the filtering distribution $\Xi_{t-1}$:
\begin{equation}
\widehat{\Xi}_{t-1}=\sum_{i=1}^{n_X} \widetilde{\omega}_{t-1}^{(i)}\delta_{X_{t-1}^{(i)}}
\end{equation}
where $\{X_{t-1}^{(i)}\}_{i=1}^{n_X} \sim \mathcal{M}\widehat{\Xi}_{t-2}$. We see from above equations that $n_X$ particles are sampled from an empirical distribution $\widehat{\Xi}_t$, i.e. it is itself defined through $n_X$ particles.

Let us give the intuition behind the reweighing step. The idea behind it is in the fact, that at this step particles with low weights have lower probability to be sampled compared with particles with high importance weights. Consequently, in this step particles with low weights will be neglected, while particles with high weights will be sampled more frequently. 

\subsection{Particle Filtering Algorithm} \label{alg0}

\hspace{3ex} 
The algorithm allows to approximate $\Xi_{t-1}^{\mu}$ by the empirical distribution $\widehat{\Xi}_{t-1}^{\mu}$ compute by the following reccurence equations:
\begin{equation}
\widehat{\Xi}_t=\widehat{\Phi}_t\widehat{\Xi}_{t-1} \ \ \mbox{with} \ \ \widehat{\Xi}_{0}=p_0
\end{equation}

where $\widehat{\Phi}_t:=\Omega_t \mathbb{M}^{n_X}\mathcal{M}$. It consists of three steps:

\begin{equation}
\widehat{\Xi}_{t-1}  \xrightarrow[]{Mutation} \mathcal{M}\widehat{\Xi}_{t-1} \xrightarrow[]{Sampling} \mathbb{M}^{n_X}\mathcal{M}\widehat{\Xi}_{t-1} \xrightarrow{Reweighing} \Omega_t \mathbb{M}^{n_X}\mathcal{M}\widehat{\Xi}_{t-1}
\end{equation}


At time $t=0$, we generate $M$ i.i.d. random variables from the prior distribution.
For $t=1,...,N-1$ we propagate $X_t \in \mathbb{R}^{n_X}$ according to the dynamics of the hidden process, update the measurement, to get a couple of random vectors $(X_{t+1},Y_{t+1})$ in the first step. Resample particles according to their probability weights $\widehat{\omega}_{t+1}(X_{t+1})$ and set resampled particles $\widehat{X}_t$.
At the final time step $t$ compute the estimate of (\ref{op1}):
\begin{equation}
\widehat{C}^{PF}=\frac{1}{n_X}\sum_{i=1}^{n_X} z(\widehat{X}_t^{(i)}) \widehat{\omega}_{t-1}^{(i)}(X_{t-1}^{(i)})\delta_{X_{t-1}^{(i)}} 
\end{equation}
where $\{X_{t-1}^{(i)}\}_{i=1}^{n_X} \sim \mathcal{M}\widehat{\Xi}_{t-2}$.

\begin{algorithm}[H] 
Initialization: $i=1,...,n_X$ - $\#$(simulations), $t=1,...,T$ - $\#$(time steps) \\
Draw $\{X_0^{(i)}\}_{i=1}^{n_X}$ from the prior $p_0(x)$. Set $\{\omega_0^{(i)}\}_{i=1}^{n_X}=\frac{1}{n_X}$;\

\For{$t=1,...,N$}{
\For{$i=1,...,n_X$}{
Propagate particles using state equation $X_t^{(i)} = f(X_{t-1}^{(i)},Y_{t-1}^{(i)},\epsilon_t)$\;
Measurement update: $Y_t = h(X_t^{(i)},Y_{t-1}^{(i)},\eta_t)$\;

Compute effective sample size $M_{eff}:\frac{1}{\sum_{i=1}^M (\omega^{(i)}_t)^2}$\;
\If{$M_{eff}<M$ or $k<N$}{
Resample using weight $\widehat{\omega}_t^{(i)}(X^{(i)}_t)\frac{\omega_t^{(i)}(X^{(i)}_t)}{\frac{1}{n_X}\sum_{j=1}^M \omega_t^{(j)}X^{(j)}_t)}$}
}
Set resampled particles as $\widehat{X}^{(i)}_t$
}

    \caption{PF Algorithm}
\end{algorithm}

Despite the advantage of sampling from highly non-linear and non-gaussian filtering distributions, we need to mention its limitations. In fact, today we have to deal with high-dimensional data, as it was shown in \cite{BLB2008}, \cite{SBB2008}, \cite{SP2015}, the collapse of weights occurs unless the sample size grows super-exponentially. Homotopy transport allows us to sample efficiently in high-dimensional framework, while avoiding the explosion of the sample size. 



\section{Homotopy Transport}

\hspace{3ex} The classical optimal transport problem is to find over all maps $\mathcal{T}:\mathbb{R}^{n_X} \rightarrow \mathbb{R}^{n_X}$, such that for $X\sim \mathbb{P}$, $\mathcal{T}(X)\sim \mathbb{Q}$ and $\mathcal{T}\in \mathcal{C}^1$; which optimizes the following criterion:
\begin{equation}
\begin{array}{c}
\inf_{\mathcal{T}} \mathbb{E}[||\mathcal{T}(X)-X||^2] \\
\mbox{s.t.} \ \ \mathbb{Q} = \mathcal{T}_{\sharp} \mathbb{P} 
\end{array}
\end{equation}

In other words, we would like to find a continuous transformation that minimizes the distance between measure $\mathbb{P}$ and measure $\mathbb{Q}$ among all these that pusheforward a prior measure $\mathbb{P}$ towards a  measure $\mathbb{Q}$. In the context of filtering problem we would like to find a transformation $\mathcal{T}$, that transport particles from a sampling measure $\mathbb{P}$ to $\mathbb{Q}$:
\begin{equation} \label{E1}
\mathbb{E}^{\mathbb{Q}}[z(X_{t})\frac{d\mathbb{P}}{d\mathbb{Q}}(X_{0:t})|\mathcal{F}_t]=\mathbb{E}^{\mathbb{P}}\left[z(\mathcal{T}(X_{t}))|\mathcal{F}_t\right]
\end{equation}

One can solve this problem using variational methods \cite{MTM2012}.
 
For the sake of exposition we represent posterior distribution, presented in the form of a normalized importance weight in the following way:
\begin{equation} \label{poster}
\psi(X_t|\textbf{Y}_{t})=\frac{1}{\mathcal{Z}_{t}}p(X_t|\textbf{Y}_{t-1})\rho(Y_t|X_t)
\end{equation}
where $\textbf{Y}_{t}=(Y_0,...,Y_t)$, the prior is $p(X_t|\textbf{Y}_{t-1})$, the likelihood is $\rho(Y_t|X_t)$ and $\mathcal{Z}_t$ is a normalization factor: $\mathcal{Z}_t=\int p(X_t|\textbf{Y}_{t-1})\rho(Y_t|X_t) dX_t$. Actually, the equation (\ref{poster}) is equivalent to the normalized importance weight in the eq. (\ref{w11}). Now, if we consider a continuous map $\mathcal{T}:\mathbb{R}^{n_X}\rightarrow \mathbb{R}^{n_X}$, then:
\begin{equation}
\psi(\mathcal{T}(X_t)|\textbf{Y}_{t})=\frac{1}{\mathcal{Z}_{t}}p(\mathcal{T}(X_t)|\textbf{Y}_{t-1})\rho(Y_t|\mathcal{T}(X_t))
\end{equation} 

Homotopy gradually modifies the prior density into the posterior density, as a scaling parameter $\lambda \in [0,1]$ increases from $0$ to $1$. In other words, by iterating we will transport homotopy $\psi(X_{t,\lambda}|Y_t)$ to a true posterior $\psi(X_t|Y_t)$, while minimizing the cost of transport. There are several conditions that homotopy has to satisfy. First, at $\lambda_0$ we should have our prior, i.e. $\psi(x_{t,\lambda_0}|Y_t)=p(X_t)$ and at some point $\lambda_{0\rightarrow 1}$, we will get approximation of our posterior density. Define a new set of density functions: $\psi(X_{t,\lambda}|\textbf{Y}_{t}):=\psi(X_t|\textbf{Y}_{t})$, $p(X_{t,\lambda}|\textbf{Y}_{t-1}):=p(X_t|\textbf{Y}_{t-1})$, $\rho(Y_t|X_{t,\lambda})^{\lambda}:=\rho(Y_t|X_{t,\lambda})$ and $\mathcal{Z}_{\lambda}:=\int p(X_{t,\lambda}|\textbf{Y}_{t-1})\rho(Y_t|X_{t,\lambda})^{\lambda} dx_{\lambda}$, so that homotopy is defined as:

\begin{equation}
\psi(X_{t,\lambda}|\textbf{Y}_{t})=\frac{1}{\mathcal{Z}_{\lambda}}\underbrace{p(X_{t,\lambda}|\textbf{Y}_{t-1})}_{prior}\underbrace{\rho(Y_t|X_{t,\lambda})^{\lambda}}_{likelihood}
\end{equation}

In order to simplify the calculation we take the logarithm of homotopy:
\begin{equation} \label{loghom}
\Psi(X_{t,\lambda}|\textbf{Y}_{t})=G(X_{t,\lambda})+\lambda L(X_{t,\lambda})-\log \mathcal{Z}_{\lambda}
\end{equation}
where $\Psi(X_{t,\lambda})=\log\psi(X_{t,\lambda}|\textbf{Y}_{t})$, $G(X_{t,\lambda}) = \log p(X_{t,\lambda}|\textbf{Y}_{t-1})$, $L(X_{t,\lambda})=\log \rho(Y_t|X_{t,\lambda})$. The dynamics of homotopy transport in the artificial time $\lambda$ is known as $log$-homotopy \cite{DH2013}. In some sense, the dynamics of transport will be given by the flow movement in the aritficial time $\lambda$, so we will look for a flow $\frac{dx}{d\lambda}$ that rules the movement of particles following log-homotopy.

If we assume that in pseudo-time $\lambda$, the flow $\frac{dx}{d\lambda}$ follows the following SDE:

\begin{equation} 
dX_{t,\lambda}=g(X_{t,\lambda})d\lambda+\eta(X_{t,\lambda})dW_{\lambda}
\end{equation}

where $W_{\lambda}$ is a vector field that pushes forward particles from prior to posterior distribution.

We impose the following assumptions:
\begin{enumerate}[I.]
\item The densities $p(X_{t,\lambda}|Y_{t-1})$ and $\rho(Y_t|X_{t,\lambda})$ are twice differentialble with respect to $X_{t,\lambda}$;
\item The function that governs the particle transport $g(X_{t,\lambda})$ is differentiable with respect to $X_{t,\lambda}$;
\item The Hessian matrix of the density $\Psi$ is non-singular;
\end{enumerate}

Now given the conditional probability density function (\ref{loghom}), we can compute the function $g(X_{t,\lambda})=\frac{dX_{t,\lambda}}{d\lambda}$ using the forward Kolmogorov equation:
\begin{equation}
\frac{\partial \psi(X_{t,\lambda})}{\partial\lambda}=-tr \left[\frac{\partial}{\partial X_{t,\lambda}}(g(X_{t,\lambda})\psi(X_{t,\lambda})) \right]+\frac{1}{2}tr\left[ \frac{\partial}{\partial X_{t,\lambda}} Q(X_{t,\lambda})\frac{\partial \psi(X_{t,\lambda})}{\partial X_{t,\lambda}}\right] 
\end{equation}

where $Q(X_{t,\lambda})=\eta(X_{t,\lambda})\eta^T(X_{t,\lambda})$ is the diffusion tensor of the process, and $tr(\cdot)$ is a trace operator. The forward Kolmogorov equation is used to relate the flow of particles $\frac{dX_{t,\lambda}}{d\lambda}$  with the evolution of  log-homotopy as $\lambda_{0\rightarrow 1}$, under the diffusion process.

\begin{multline} \label{1}
\frac{\partial \psi(X_{t,\lambda})}{\partial\lambda}=-tr \left[\psi(X_{t,\lambda})\frac{\partial g(X_{t,\lambda})}{\partial X_{t,\lambda}}+g(X_{t,\lambda})^T\frac{\partial \psi(X_{t,\lambda})}{\partial X_{t,\lambda}} \right]+\frac{1}{2}div\left[ \frac{\partial}{\partial X_{t,\lambda}} Q(X_{t,\lambda})\frac{\partial \psi(X_{t,\lambda})}{\lambda}\right] = \\ =  -\psi(X_{t,\lambda})tr\left[\frac{\partial g(X_{t,\lambda})}{\partial X_{t,\lambda}}\right]-g(X_{t,\lambda})^T\frac{\partial \psi(X_{t,\lambda})}{\partial X_{t,\lambda}}+\frac{1}{2}div\left[ \frac{\partial}{\partial X_{t,\lambda}} Q(X_{t,\lambda})\frac{\partial \psi(X_{t,\lambda})}{\partial X_{t,\lambda}}\right]
\end{multline}
where $div(\cdot)$ is a divergence operator.
On the other hand if we take the derivative of equation (\ref{loghom}) with respect to $\lambda$, we have: 
\begin{equation} \label{11}
\frac{\partial \Psi(X_{t,\lambda})}{\partial \lambda}=L(X_{t,\lambda})-\frac{\partial}{\partial \lambda}\log \mathcal{Z}_{\lambda}
\end{equation}

Since $\Psi(X_{t,\lambda})$ is a composition of two functions, we will need to use the chain rule:
\begin{equation}  \label{12}
\frac{\partial \Psi(X_{t,\lambda})}{\partial \lambda}=\frac{1}{\psi(X_{t,\lambda})}\frac{\partial \psi(X_{t,\lambda})}{\partial \lambda}
\end{equation}

By substituting eq. \eqref{12} into \eqref{11} and rearranging the terms:
\begin{equation} \label{2}
\frac{\partial \psi(X_{t,\lambda})}{\partial \lambda}=\psi(X_{t,\lambda})\left[L(X_{t,\lambda})-\frac{\partial}{\partial \lambda}\log \mathcal{Z}_{\lambda}\right]
\end{equation}

Observe that \eqref{1} and \eqref{2} are identical, so by equating and dividing on $\psi(X_{t,\lambda})$ we get:
\begin{multline} \label{13}
L(X_{t,\lambda})-\frac{\partial}{\partial \lambda}\log \mathcal{Z}_{\lambda}=-g(X_{t,\lambda})^T\frac{1}{\psi(X_{t,\lambda})}\frac{\partial \psi(X_{t,\lambda})}{\partial X_{t,\lambda}}-\\- 
tr\left[\frac{\partial g(X_{t,\lambda})}{\partial X_{t,\lambda}} \right]+\frac{1}{2\psi(X_{t,\lambda})}div\left[ \frac{\partial}{\partial X_{t,\lambda}} Q(X_{t,\lambda})\frac{\partial \psi(X_{t,\lambda})}{\partial X_{t,\lambda}}\right]
\end{multline}

In \cite{DH2015}, authors propose to take the derivative of \eqref{13} with respect to $X_{t,\lambda}$ in order to find explicitely the equation of flow on the one hand, and to get rid of the normalization constant $\mathcal{Z}_{\lambda}$ that lead to instabilities on the other hand.

\begin{multline} \label{14}
\frac{\partial L(X_{t,\lambda})}{\partial X_{t,\lambda}}=-g(X_{t,\lambda})^T \frac{\partial^2 \Psi(X_{t,\lambda})}{\partial X_{t,\lambda}^2}-\frac{\partial \Psi(X_{t,\lambda})}{\partial X_{t,\lambda}}\frac{\partial g(X_{t,\lambda})}{\partial X_{t,\lambda}}-\frac{\partial}{\partial X_{t,\lambda}}tr\left[\frac{\partial g(X_{t,\lambda})}{\partial X_{t,\lambda}} \right]+\\+\frac{\partial }{\partial X_{t,\lambda}}\left(\frac{1}{2\psi(X_{t,\lambda})}div\left[Q(X_{t,\lambda})\frac{\partial \psi(X_{t,\lambda})}{\partial X_{t,\lambda}}\right]\right)
\end{multline}

Observe that we get a highly nonlinear PDE. We use the solution found in \cite{DH2013} and \cite{DH2011}, which states that if we could find a vector field $g(X_{t,\lambda})$ and diffusion tensor $Q(X_{t,\lambda})$, such that sum of the last three terms in \eqref{14} are equal to zero. The PDE, then simplifies to:

\begin{equation}
\frac{\partial L(X_{t,\lambda})}{\partial X_{t,\lambda}}=-g(X_{t,\lambda})^T \frac{\partial^2 \Psi(X_{t,\lambda})}{\partial X_{t,\lambda}^2}
\end{equation}

Using the assumption III, i.e. the Hessian matrix $\frac{\partial^2 \Psi(X_{t,\lambda})}{\partial X_{t,\lambda}^2}$ is non-singular, we get explicitely the flow $g(X_{t,\lambda})$:
\begin{equation}
g(X_{t,\lambda})=-\left[\frac{\partial^2 \Psi(X_{t,\lambda})}{\partial X_{t,\lambda}^2}\right]^{-1} \left[\frac{\partial L(X_{t,\lambda})}{\partial X_{t,\lambda}} \right]^T
\end{equation} 

\subsection{Homotopy Transport Algorithm} \label{alg1}

\hspace{3ex}  \textbf{Sampling from the prior.}
First we generate $M$ i.i.d random variables $X^{(i)}_t$ from the prior density $p_0(x)$, initialize pseudo-time $\lambda$ and  set the state variables that will be transported as: $X_{t,\lambda}^{(i)}=X_{t|t-1}^{(i)}$.

\textbf{Transportation Stage.}
For $t=2,...,N-1$, compute the derivative with respect to $X_{t,\lambda}$ of the measurement function. If $h$ is non-linear, a second order Taylor expansion at $X_{t,\lambda}$ allows speeding up the calculation by linearizing the first derivative. After that, update the pseudo time by setting : $\lambda=\lambda + \Delta \lambda$.

Compute the flow $g(X_{t,\lambda}^{(i)})$. Note, that the first Hessian could be derived by twice differentiating a log-homotopy equation (\ref{loghom}):
\begin{equation} \label{Psi}
\frac{\partial^2 \Psi(X_{t,\lambda}^{(i)})}{\partial X_{t,\lambda}^2}=\frac{\partial^2 G(X_{t,\lambda}^{(i)})}{\partial X_{t,\lambda}^2}+\lambda \frac{\partial^2 L(X_{t,\lambda}^{(i)})}{\partial X_{t,\lambda}^2}
\end{equation} 

The first term in (\ref{Psi}) $\frac{\partial^2 G(X_{t,\lambda}^{(i)})}{\partial X_{t,\lambda}^2}$ is estimated by using a sample covariance matrix of $t$ patricles generated form the prior distribution:
\begin{equation}
\frac{\partial^2 G(X_{t,\lambda})}{\partial X_{t,\lambda}^2} \approx - \widehat{S}^{-1}_{M_x}
\end{equation}
Compute the transportation of particles from the measure $\mathbb{P}$ to the measure $\mathbb{Q}$:
\begin{equation}
X^{(i)}_{t,\lambda}=X^{(i)}_{t,\lambda}+\Delta \lambda g(X_{t,\lambda}^{(i)})
\end{equation} 
And finally update the state parameter: 
\begin{equation}
\breve{X}_t=\frac{1}{n_X}\sum_{i=1}^{n_X} X^{(i)}_{t,\lambda}
\end{equation} 

\textbf{Maturity.}

At the final time interval $]N-1,N]$ compute the estimator of (\ref{E1}):
\begin{equation}
\widehat{C}^{HT}=\frac{1}{n_X}\sum_{i=1}^{n_X} z(\breve{X}_t^{(i)})
\end{equation}

\begin{algorithm}[H] 
Initialization: $i=1,...,n_X$ - $\#$(simulations), $t=1,...,N$ - $\#$(time steps) \\
Draw $\{X_0^{(i)}\}_{i=1}^{n_X}$ from the prior $p_0(x)$.\\
Set $\{\omega_0^{(i)}\}_{i=1}^{n_X}=\frac{1}{n_X}$

\For{$t=1,...,N$}{
\For{$i=1,...,n_X$}{
Propagate particles using state equation $X_t^{(i)} = f(X_{t-1}^{(i)},Y_{t-1}^{(i)},\epsilon_t)$\;
Measurement update: $Y_t = h(X_t^{(i)},Y_{t-1}^{(i)},\eta_t)$\;
Initialize pseudo-time $\lambda=0$\;
Set $X_{t,\lambda}^{(i)}=X_{t|n-1}^{(i)}$\;

\While{$\lambda < 1$}{

Compute SCM $\widehat{S}_M$;\

Calculate an estimate: $X_{t,\lambda}=\frac{1}{n_X}\sum_i X_{t,\lambda}^{(i)}$

Compute the matrix $\widehat{H}=\frac{\partial h(X_{t,\lambda}^{(i)})}{\partial X_{t,\lambda}}$\;

Update the time: $\lambda=\lambda+\Delta \lambda$\;
Calculate the flow $\frac{dX^{(i)}_{t,\lambda}}{d\lambda}=-\left[\frac{\partial^2 \Psi(X^{(i)}_{t,\lambda})}{\partial X_{t,\lambda}^2}\right]^{-1} \left[\frac{\partial L(X^{(i)}_{t,\lambda})}{\partial X_{t,\lambda}} \right]^T$\;
Transport particles according to its flow: $X_{t,\lambda}^{(i)}=X_{t,\lambda}^{(i)}+\Delta \lambda \frac{d X_{t,\lambda}^{(i)}}{d\lambda}$;\
}

Update state estimate:\\
$\breve{X}_t=\frac{1}{n_X}\sum_{i=1}^{n_X} X_{t,\lambda}^{(i)}$

}
}

    \caption{Homotopy Transport Algorithm}
\end{algorithm}

\section{Homotopy Transport with Particle Reweighing}

\hspace{3ex}  Taking into account the difficulties one faces in non-Gaussian and high-dimensional problems, the idea of a particle transport without any use of sampling techniques is very elucidating. The next question that arises is whether we could direct the transportation by choosing those particles that have higher probability of reaching rarely visited areas of the state space? We propose a mixture of homotopy particle transport with a particle reweighing at each time step. The numerical test that we performed on the toy example of a Stein-Stein stochastic volatility model showes that we significantly reduce the variance and bias of our estimator.   

The algorithm consists of two steps: first we transport particles according to its flow, and second, we choose those particles that have higher probability of faster exploring the state space. 

\begin{equation} \label{RN4.1}
\mathbb{E}^{\widetilde{\mathbb{Q}}}[z(X_{t})\frac{d\mathbb{P}}{d\widetilde{\mathbb{Q}}}(X_{0:t})|\mathcal{F}_t]=\mathbb{E}^{\mathbb{P}}\left[z(\mathcal{T}(X_{t}))|\mathcal{F}_t\right]
=\mathbb{E}^{\mathbb{Q}}\left[z(\mathcal{T}(X_{t}))\frac{d\mathbb{P}}{d\mathbb{Q}}(X_{0:t})| \mathcal{F}_t\right]
\end{equation}

where $\mathcal{T}$ is a flow of particles under the pseudotime $\lambda$ decribed in the section $\ref{alg1}$.

By setting $\textbf{X}_t=(X_0,...,X_t)$, we could express our Radon-Nikodym derivative in a product form:
\begin{equation}
\frac{d\mathbb{P}}{d\mathbb{Q}}(X_{0:t})=\frac{d\mathbb{P}}{d\widetilde{\mathbb{Q}}}\times \frac{d\widetilde{\mathbb{Q}}}{d\mathbb{Q}}(X_{0:t})
\end{equation} 

where the first Radon-Nikodym derivative denotes the transport of particles from a mesure $\mathbb{P}$ to a measure $\widetilde{\mathbb{Q}}$, then we choose the particles that have high probability of reaching rare corners of the state space, using $\frac{d\widetilde{\mathbb{Q}}}{d\mathbb{Q}}$ that allows us to reassess the weights of the particles.


As in the section 2, an importance measure $\mathbb{Q}$ that will play a resampling to choose the trajectories with higher weight, given that there exists an importance kernel $\widetilde{K}_t$, such that $K_t\ll \widetilde{K}_t$, could be defined as:
\begin{equation}
\mathbb{Q}(B)=\sum_{i=1}^{n_X} \omega^{(i)}_t \widetilde{K}_t(X^{(i)}_t,A_i), \ \ B\in\mathcal{B}(\mathbb{R}^{n_X})
\end{equation}

where the set $A_i=\{\mathcal{T}(X_{t+1})\in\mathbb{R}^{n_X}|\mathbb{1}_{B}(X^{(i)}_t,\mathcal{T}(X_{t+1}))=1\}$. Assuming, that the prior and sampling kernels $K_t$ and $\widetilde{K}_t$ have densities $k_t$ and $\widetilde{k}_t$ respectively, then the Radon-Nikodym derivative is
\begin{equation} \label{RN}
\frac{d\widetilde{\mathbb{Q}}}{d\mathbb{Q}}(X_{0:t})=\prod_{l=0}^{t} \rho_{l}(\mathcal{T}(X_{l}),Y_{l})\frac{ \omega_{l-1}(X_{l-1})k_l(X_{l-1},\mathcal{T}(X_{l}))}{\omega_{l-1}(X_{l-1})\widetilde{k}_l(X_{l-1},\mathcal{T}(X_{l}))}
\end{equation}
such that $\omega_t(X_t)=\omega_t^{(i)}(X_t^{(i)})$ if $X_t=X_t^{(i)}$, and $\omega_t(X_t)=1$ otherwise.

The an unnormalized weight is given by:

\begin{equation}
\omega_{t}^{(i)}(\mathcal{T}(X_{t}^{(i)}))=\prod_{l=1}^t\rho_{l}(\mathcal{T}(X_{l}^{(i)}),Y_{l})\frac{ k_l(X_{l-1},\mathcal{T}(X_{l}^{(i)}))}{\widetilde{k}_l(X_{l-1},\mathcal{T}(X_{l}^{(i)}))}
\end{equation} 

So, now we have homotopy transport with particle reweighing estimator:

\begin{equation}
\widehat{C}^{TRW}=\frac{1}{n_X}\sum_{i=1}^{n_X} z(\mathcal{T}(X_t^{(i)})) \widehat{\omega}_{t-1}^{(i)}(\mathcal{T}(X_{t-1}^{(i)}))
\end{equation}

\subsection{PF-Enhanced Homotopy Transport Algorithm}

\hspace{3ex}  The  algorithm could be described by the following scheme, $\forall i = 1,...,n_X$:
\begin{equation}
X_t^{(i)}  \xrightarrow[]{Sampling} X_{t+1}^{(i)} \xrightarrow[]{Transportation} \mathcal{T}(X_{t+1}^{(i)})=\breve{X}_{t+1}^{(i)} \xrightarrow{Reweighing} \Phi(\breve{X}_{t+1}^{(i)})= \widehat{X}_{t+1}^{(i)}
\end{equation}

where $\Phi$ is an operator that denotes the resampling mechanism of particles. 
If we assume that there is a continuous kernel $\widetilde{K}_t$, such that $K_t\ll \widetilde{K}_t$ with densities $k_t$ and $\widetilde{k}_t$ respectively, then we can define a weight function $\omega_{t}^{(i)}$:
\begin{equation}
\omega_{t}^{(i)}(\breve{X}_{t}^{(i)})=\prod_{l=1}^t\rho_{l}(\breve{X}_{l}^{(i)},Y_{l})\frac{ k_l(\widehat{X}_{l-1}^{(i)},\breve{X}_{l}^{(i)})}{\widetilde{k}_l(\widehat{X}_{l-1}^{(i)},\breve{X}_{l}^{(i)})}
\end{equation}

\subsubsection{Detailed Algorithm}

\hspace{3ex} \textbf{Sampling from the prior.}
As in the section \ref{alg1}, we start with $M$ particles sampled from the prior distribution $p_0$, initialize pseudo-time $\lambda$ and  set the state variables that will be transported as: $X_{t,\lambda}^{(i)}=X_{t|n-1}^{(i)}$.

\textbf{Transportation Stage.}
Follow steps 6-8 of the Algorithm 2 in the section \ref{alg1}.

\textbf{Path Reweighing Stage.}
Compute the normalized importance weight:

\begin{equation}
\widehat{\omega}_{t}^{(i)}(\breve{X}_{t}^{(i)})=\frac{\omega_{t}^{(i)}(\breve{X}^{(i)}_{t})}{\sum_{i=1}^{n_X} \omega_{t}^{(j)}(\breve{X}_{t}^{(j)})}
\end{equation}

\textbf{Maturity}
At the time interval $]N-1,N]$ compute the final Homotopy transport reweighted estimator:
\begin{equation}
\widehat{C}^{TRW}=\frac{1}{n_X}\sum_{i=1}^{n_X} z(\widehat{X}_t^{(i)}) \widehat{\omega}_{t-1}^{(i)}(\breve{X}_{t-1}^{(i)})
\end{equation}

\begin{algorithm}[H] 
Initialization: $i=1,...,n_X$ - $\#$(simulations), $t=1,...,T$ - $\#$(time steps) \\
Draw $\{X_0^{(i)}\}_{i=1}^{n_X}$ from the prior $p_0(x)$.\\
Set $\{\omega_0^{(i)}\}_{i=1}^{n_X}=\frac{1}{n_X}$

\For{$t=1,...,N$}{
\For{$i=1,...,n_X$}{
Follow steps 6-8 of the Algorithm 2 in the section \ref{alg1}.

Follow steos 7-12 of the Algorithm 1 in the section \ref{alg0}
}
}

    \caption{Homotopy Transport with Particle Reweighing Algorithm}
\end{algorithm}

\section{Numerical Applications and Results}
\hspace{3ex}  As a toy example, we decided to test the algorithms on a Stein-Stein stochastic volatility model. We set log-returns as $Y_t=\log(S_t)$, then the model takes the following form:
\begin{equation}
\left\{\begin{array}{c} 
dY_t=(\mu-\frac{X_t^2}{2})dt+X_tdB_t\\
dX_t=\kappa(\theta-X_t)dt+\sigma dW_t

\end{array} \right.
\end{equation}

where $X_t$ is a volatility process, $Y_t$ the dynamics of log-returns, $\mu$ is a drift, $\theta$ is a long-term variance, $\kappa$ - the rate of reversion, $\sigma$ is the volatility of volatility, and $B_t$ and $W_t$ are two independent Brownian motions, in the sense that $\langle dB_t,dW_t\rangle=0$. 

Using the above presented stochastic volatility model, we would like to compute estimates for a European option. For a given interest rate $r$, maturity $T$, strike price $K$, and for a function $z(\cdot,x)=\max(x-K,0)$, the call price of the option is given by:
\begin{equation} 
C(X_t,Y_t)=B_{t,T}\mathbb{E}^{\mathbb{P}}\left[z(X_T,Y_T)| \mathcal{F}_t\right]
\end{equation}
where $\mathcal{F}_t=\sigma\{(Y_0,...,Y_t)\}$.

We chose Euler-Muruyama discretization scheme, which gives:

\begin{equation}
\left\{\begin{array}{c} 
Y_t-Y_{t-1}=(\mu-\frac{X_{t-1}^2}{2})\Delta t+X_{t-1} \sqrt{\Delta t} \epsilon_t\\
X_t-X_{t-1}=\kappa(\theta-X_{t-1})\Delta t+\sigma  \sqrt{\Delta t} \eta_t

\end{array} \right.
\end{equation}
where $\Delta t$ is a discretization size, $\epsilon_t$ and $\eta_t$ are independent Gaussian variates, $\mathcal{N}(0,1)$.

We compare each approach by estimating the standard deviations, the root mean squared error (RMSE), the bias, the relative mean squared error(RRMSE), the time required to compute each estimate and the figure of merit (FOM). We run 20 Monte Carlo experiments. For $l=1,...,M_s$ the RMSE estimator is given by:
\begin{equation}
RMSE = \sqrt{\frac{1}{M_s}\sum_{l=1}^{M_s} ||C-\widehat{C}_{l}||^2}
\end{equation}
where $C$ is the price computed analytically, $\widehat{C}_{l}$ are Monte Carlo estimates and $M_s=20$ is the number of Monte Carlo experiments. As a reference price, we used the article by EM Stein \cite{SS1991}.

\begin{equation}
Bias = \sqrt{RMSE^2 - St.dev^2}
\end{equation}

where $St.dev$ are standard deviations of MC estimates. The RRMSE is computed using the following formula:
\begin{equation}
RRMSE = \frac{RMSE}{\widehat{C}}
\end{equation}

To measure the efficiency of each method presented in the article, we will use the figure of merit(FOM)\cite{RT2009}:
\begin{equation}
FOM = \frac{1}{R^2\times CPU_t}
\end{equation}

where $CPU_t$ is the CPU time need to compute the estimator and $R$ is a relative error, which is the measure of statistical precision:
\begin{equation}
R = \frac{St. dev}{\bar{C}} \propto \frac{1}{\sqrt{M}}
\end{equation}
where $\bar{C}=\sum_{l=1}^{M_s} \widehat{C}_{l}$

We used 20 000 and 40 000 simulations over 64 time intervals for our MC experiments. Table 1. shows that homotopy and reweighted(RW)-homotopy algorithms shows less statistical errors then traditional PF. If we compare homotopy and RW-homotopy, we could see that FOM  says that the first is more efficient the latest, due to the fact that we need more time to reweight the paths. Meanwhile RW-homotopy shows less erros and st. deviations.

\begin{table}[H] 
 \caption{Stein-Stein Stochastic volatility option price estimates statistics. $S_0=100, \ K=90, \ r=0.0953, \sigma = 0.2, \ \kappa=4, \ \theta=0.25, \ V_0=0.25, \ T=1/2$, and dividends $d=0$ True price: $16.05$, $t=20000$, $M=64$}                                    

\centering                                                                          
 \begin{tabular}{|c|c|c|c|c|}                                                     
\hline                                                                       
Stat  & MC & PF &   Homotopy &    RW-Homotopy   \\
 
\hline
St. dev. & 0.127495344 & 0.106264197 & 0.102775848 & 0.08360908\\
\hline                                                                      
RMSE &  0.148073563 & 0.115032508 & 0.105302932 & 0.084510606\\
\hline                                                                      
Bias & 0.075304165 & 0.044049955 & 0.022931037 & 0.012311146\\
\hline
RRMSE &  0.00137298 & 0.000827032 & 0.000827032  & 0.000444367 \\
\hline
CPU time &  0.1327525 & 0.31177 &  0.179 & 0.38819 \\
\hline
FOM &  118181.69 & 72715.84 & 135692.61 & 95193.97 \\
\hline
\end{tabular}
\end{table}

\begin{figure}[H]
		\centering
	\includegraphics[width=0.7\linewidth]{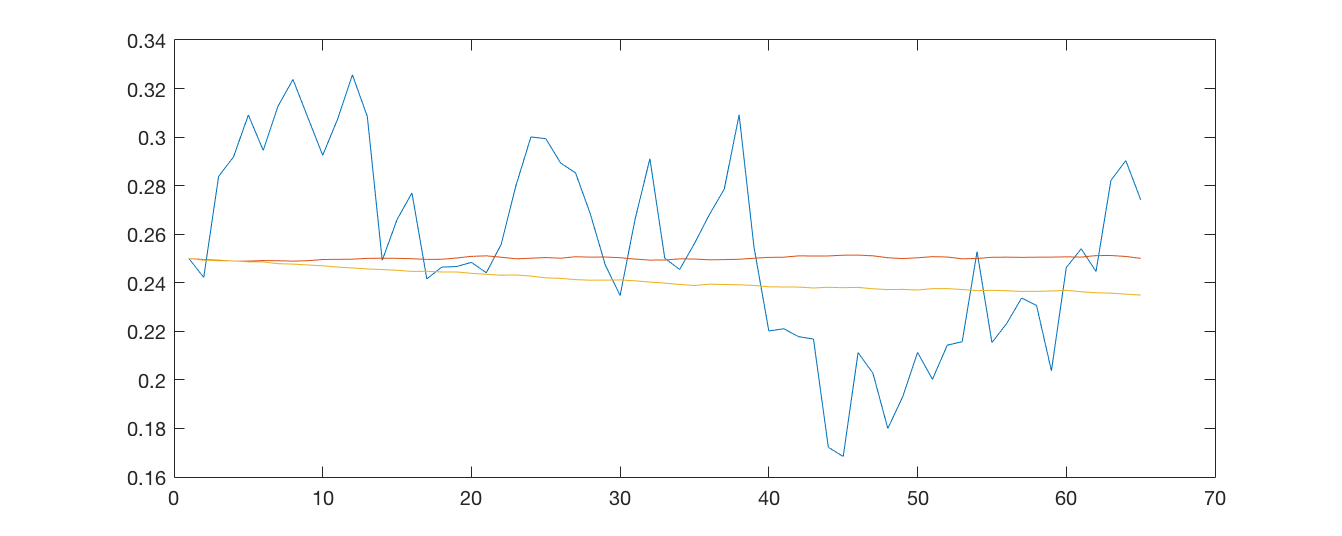}	
	\caption{Volatily dynamics, PF (Blue), Homotopy (Red), RW-homotopy(Yellow)}
	\end{figure}
	
\begin{figure}[H]
	\includegraphics[width=0.49\linewidth]{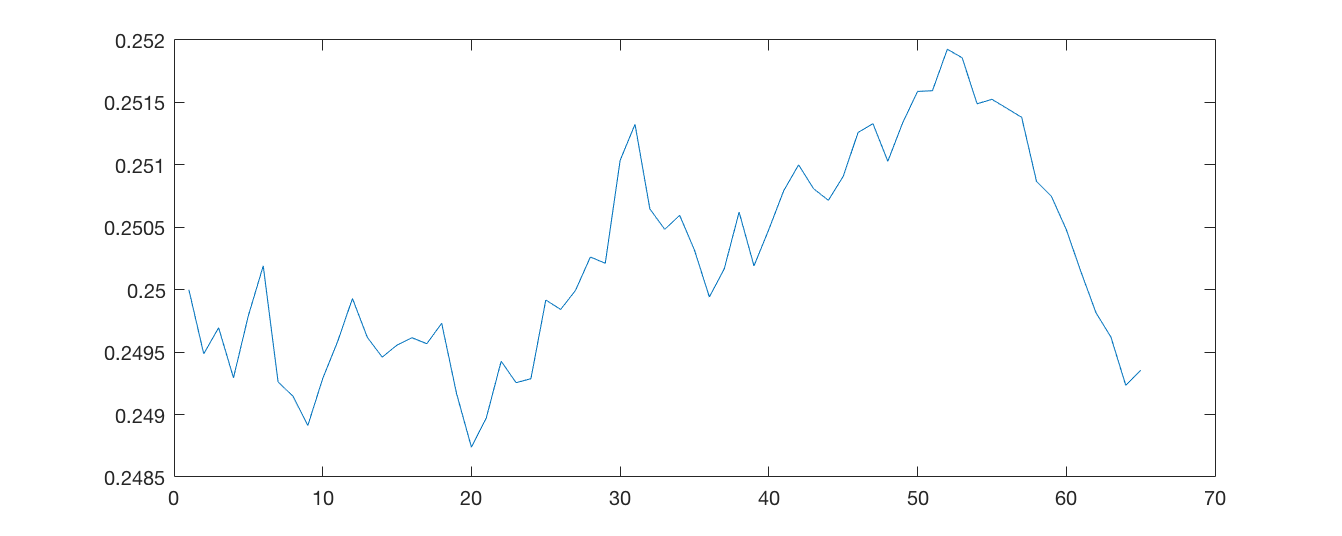}
	\includegraphics[width=0.49\linewidth]{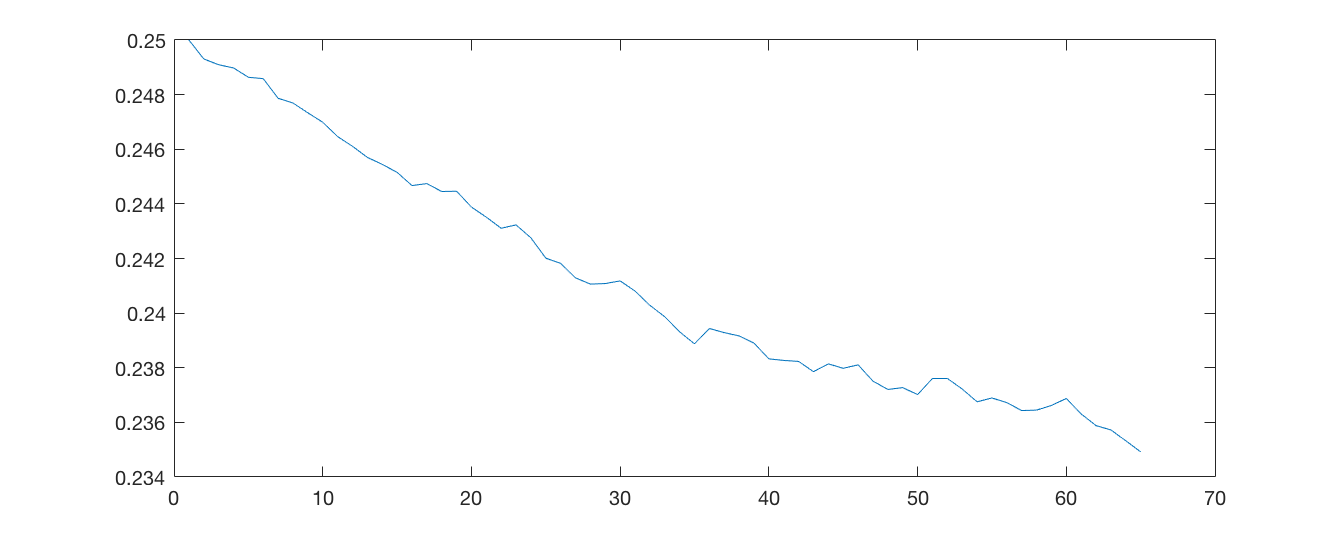}
		\caption{Zoomed volatilty dynamics. Homotopy (left), RW-homotopy (right)}
	\end{figure}

\begin{table}[H]
 \caption{Stein-Stein Stochastic volatility option price estimates statistics. $S_0=100, \ K=90, \ r=0.0953, \sigma = 0.2, \ \kappa=4, \ \theta=0.25, \ V_0=0.25, \ T=1/2$, and dividends $d=0$ True price: $16.05$, $t=40000$, $M=64$}                                    

\centering                                                                          
 \begin{tabular}{|c|c|c|c|c|}                                                     
\hline                                                                       
Stat  & MC & PF &   Homotopy &    RW-Homotopy   \\
 
\hline
St. dev. & 0.070351719 & 0.060799052 & 0.048943672 & 0.045246118\\
\hline                                                                      
RMSE & 0.130446299 & 0.079273246 & 0.04921257 & 0.045762201 \\
\hline                                                                      
Bias & 0.109849318 &0.050869665 & 0.005137504 & 0.006853309\\
\hline
RRMSE & 0.001067583 & 0.000392831 & 0.00015101 & 0.000130578\\
\hline
CPU time & 0.278895 & 0.54737 & 0.26618 & 0.581495\\
\hline
FOM & 184049.069 & 126479.8136  & 403391.758 & 216062.7397\\
\hline
\end{tabular}
\end{table}

Despite the fact that Monte Carlo estimate showed higher FOM, than PF, due to the fact that it takes less time to compute Monte Carlo estimator. Whereas PF has lower RMSE and the bias.

\section{Conclusions and Further Research}

\hspace{3ex}  The estimation of latent variables has a lot of applications in engineering and finance. We provide homotopy based algorithm and its extension with reweighted trajectories that permits to solve the optimal transportation problem. 

Numerical results that we applied in European option pricing with stochastic volatility demonstrated the efficiency of  the proposed algorithms with respect to error, bias and other stastics. Both algorithms ourperformed Particle filtering. The path-reweighing allowed to reduce standard deviations, and in some cases the bias and the RMSE compared to the homotopy transport algorithm.

From our experiments we could observe the following:
\begin{itemize}
\item Homotopy transport is fast algorithm, which is spectacularily demonstrated in the figure of merit statistics.
\item Efficiency of homotopy transport algorithm increases as the number of particles increases.
\item Implementation of homotopy transport requires less effort than a vanilla Monte Carlo algorithm.
\item Homotopy transport proved to be unbiased estimator.
\item Homotopy with path reweighing proved to reduce the bias when the number of particles is small compared to homotopy transport without reweighing. 
\end{itemize}

While reweighted homotopy transport approach showed the reduced RMSE and Bias in low-dimensions, the mixture of homotopy transport and bootstrap resampling, the importance weight could converge to unity in high-dimensional problems(\cite{SP2015}). In the next article we plan to investigate this issue. It will be also interesting to check the homotopy transport on non-gaussian examples.

\section*{Appendix}
\subsection*{Flow related computations}

In the classical particle filtering approach the  desired estimate is approximate by $M$ pratciles:
\begin{equation}
p(X_t|Y_{1:n-1})=\frac{1}{n_X}\sum_{i=1}^{n_X} p(X_t|X_{t-1}=X_{t-1}^{(i)})
\end{equation}
So that posterior at time $t$:
\begin{equation}
p(X_t|Y_{1:n})=\frac{1}{\mathcal{Z}_t}
\rho(Y_t|X_t)p(X_t|Y_{1:n-1})
\end{equation}
The transition density is given by:
\begin{equation}
p(X_t|X_{t-1})=\mathcal{N}(X_t;\mu_X,\sigma_X)
\end{equation}
where $\mu_X=X_{t-1}+\kappa(\mu_X-X_{t-1})\Delta t$ and $\sigma_X=\sigma_X^2 \Delta t$.

The likelihood, $p(Y_t|Y_{t-1},X_{t-1})$:
\begin{equation}
m_{t,N}(X_{t-1})=Y_{t-1}+(\mu-\frac{X_{t-1}^2}{2})\Delta t
\end{equation}
\begin{equation}
\sigma^p_{t,N}(X_{t-1})=X_{t-1}^2\Delta t
\end{equation}

So,
\begin{equation}
\rho(Y_t|Y_{t-1},X_{t-1})=\prod_{t=1}^N\mathcal{N}(Y_t;m_{t,N}(X_{t-1}),\sigma^p_{t,N}(X_{t-1}))
\end{equation}
The unnormalized posterior is given by:
\begin{equation}
\mathcal{P}_t=\rho(Y_t|Y_{t-1},X_{t-1})
\end{equation}
\begin{equation}
X_t=\psi(X_t;\mathcal{P}_t)
\end{equation}

Next,
\begin{equation}
\psi(X)=-\log(\mathcal{P}(X))
\end{equation}

by removing some constants that have no impact on posterior distribution, we have
\begin{equation}
\psi(X)=\sum_{t=1}^N \frac{(Y_t-m_{t,N}(X_{t-1}))^2}{2\sigma^p_{t,N}(X_{t-1})}+\frac{1}{2}\log(\sigma^p_{t,N}(X_{t-1}))
\end{equation}

\begin{multline}
\frac{\partial \psi}{\partial x}(X)=\frac{1}{2} \sum_{t=1}^N\left( \frac{\nabla_X \sigma^p_{t,N}(X_{t-1})}{X_{t-1}} \right. - \\ - \left. \frac{(Y_t-m_{t,N})(2\sigma^p_{t,N}(X_{t-1})\nabla_X m_{t,N}(X_{t-1})+(Y_t-m_{t,N})\nabla_X \sigma^p_{t,N}(X_{t-1}))}{\sigma^p_{t,N}(X_{t-1})^2}  \right)
\end{multline}

\begin{equation}
\nabla_X m_{t,N}(X_{t-1})=-X_{t-1}\Delta t
\end{equation}
\begin{equation}
\nabla_X \sigma^p_{t,N}(X_{t-1}))=2X_{t-1}\Delta t
\end{equation}

\begin{multline}
\frac{\partial \psi}{\partial x}(X)=\frac{1}{2} \sum_{t=1}^N\left( 2\Delta t - \right. \\ - \left. \frac{(Y_t-m_{t,N})(2\sigma^p_{t,N}(X_{t-1})\nabla_X m_{t,N}(X_{t-1})+(Y_t-m_{t,N})\nabla_X \sigma^p_{t,N}(X_{t-1}))}{\sigma^p_{t,N}(X_{t-1})^2}  \right)
\end{multline}

\begin{equation}
u=(Y_t-m_{t,N})(-2X_{t-1}^3\Delta t^2+2(Y_t-m_{t,N})X_{t-1}\Delta t)
\end{equation}
\begin{equation}
u'=2X_{t-1}^2\Delta t^2 \left(X_{t-1}^2 \Delta t - (Y_t-m_{t,N})\right)+(Y_t-m_{t,N})\left(-6X^2_{t-1}\Delta t^2 + 2\Delta t((Y_t-m_{t,N}))\right)
\end{equation}
\begin{equation}
v= X^4_{t-1}\Delta t^2,\ \ \ v'=4X_{t-1}^3\Delta t^2
\end{equation}
\begin{equation}
\frac{\partial^2 \psi}{\partial x^2}=\frac{u'v-v'u}{v^2}
\end{equation}


\begin{thebibliography} {1}
\bibitem{BLB2008} P. Bickel, B. Li, and T. Bengtsson, “Sharp failure rates for the bootstrap particle filter in high dimensions,” Institute of Mathematical Statistics
Collections, vol. 3, pp. 318–329, 2008.
\bibitem{DH2013} Daum, F., \& Huang, J. (2013). Particle flow with non-zero diffusion for nonlinear filters. In Proceedings of spie conference (Vol. 8745).
\bibitem{DH2011} Daum, F., \& Huang, J. (2011). Particle degeneracy: root cause and solution. In Proceedings of spie conference (Vol. 8050).
\bibitem{DH2015} Daum, F., \& Huang, J. (2015). Renormalization group flow in k-space for nonlinear filters, Bayesian decisions and transport.
\bibitem{DM2013} Del Moral, P.: Mean field simulation for Monte Carlo integration. CRC Press (2013)
\bibitem{DM2004}  Del Moral, P.: Feynman-Kac Formulae: Genealogical and Interacting Particle Systems with Applications. Probability and Applications. Springer, New York (2004).
\bibitem{DM} Del Moral, P.: Nonlinear Filtering: Interacting Particle Solution(1996). Markov Processes and Related Fields 2 (4), 555-580
\bibitem{MTM2012} El Moselhy, Tarek A. and Marzouk, Youssef M.(2012). Bayesian inference with optimal maps. Journal of Computational Physics. (Vol. 231)
\bibitem{SS1991} Stein, Elias M, and Jeremy C Stein. 1991. “Stock Price Distributions with Stochastic Volatility: An Analytic Approach.” Review of Financial Studies 4: 727-752.
\bibitem{RT2009} Rubino, G., Tuffin, B.: Rare event simulation using Monte Carlo methods. Wiley (2009)
\bibitem{HPL2013} Beiglböck, M., Henry-Labordère, P. \& Penkner, F. Finance Stoch (2013) 17: 477. doi:10.1007 /s00780-013-0205-8
\bibitem{CV} Villani, C. : Topics in optimal transportation, Graduate studies in Mathematics AMS, Vol 58.
\bibitem{Rachev} Rachev, S. T. and Ruschendorf, L. : Mass Transportation Problems. In Vol. 1: Theory. Vol. 2: Applications. Springer, Berlin, 1998.
\bibitem{SBB2008} C. Snyder, T. Bengtsson, P. Bickel, and J. Anderson, “Obstacles to high-dimensional particle filtering,” Monthly Weather Review, vol. 136, no. 12,
pp. 4629–4640, 2008.
\bibitem{SP2015} F. Septier and G. W. Peters, “An Overview of Recent Advances in Monte-Carlo Methods for Bayesian Fitlering in High-Dimensional Spaces,”
in Theoretical Aspects of Spatial-Temporal Modeling, G. W. Peters and T. Matsui, Eds. SpringerBriefs - JSS Research Series in Statistics,
2015.
\end{thebibliography}
\end{document}